\documentclass{amsart}

\usepackage{graphicx}

\usepackage[mode=buildnew]{standalone}
\usepackage{float}

\theoremstyle{definition}

\theoremstyle{remark}

\numberwithin{equation}{section}

\begin{document}

\title[Automated Math and the Reconfiguration of Proof and Labor]{Automated Mathematics and the Reconfiguration of Proof and Labor}

\author{Rodrigo Ochigame}
\address{Institute of Cultural Anthropology and Development Sociology, Leiden University, Leiden, The Netherlands}
\curraddr{}
\email{rodrigo@ochigame.org}
\thanks{}



\dedicatory{}

\begin{abstract}
This essay examines how automation has reconfigured mathematical proof and labor, and what might happen in the future. It discusses practical standards of proof, distinguishes between prominent forms of automation in research, provides critiques of recurring assumptions, and asks how automation might reshape economies of labor and credit.
\end{abstract}

\maketitle

\section{Introduction}

Dreams of automating mathematical research are as old as the imagination of early mechanical calculating devices, long before the creation of modern digital computers \cite{daston_enlightenment_1994,jones_reckoning_2016,ames_call_2022}. The latter have already generated diverse aspirations to automated mathematics, from formal verification to interactive and automated theorem proving, since the mid-twentieth century. Yet, novel developments in the twenty-first are pressing mathematicians to reconsider such ideas and their implications afresh. Formal verification, once imagined to work only for rudimentary results, now seems feasible even for some of the most sophisticated results of contemporary research. Advances in artificial intelligence and machine learning, which have captured massive amounts of public attention and capital investment, promise to expand the automation of theorem proving.

In what follows, I examine how automation has reconfigured mathematical proof and labor, and what might happen in the future. My perspective is grounded empirically in the comparative analysis of a wide range of historical and contemporary cases. First, I suggest that past and present controversies about the status of computer-assisted proofs reflect a longstanding tension in modern mathematics, between what I call post-Cartesian and post-Leibnizian practical standards of proof. Then, I distinguish between prominent forms of automation in mathematical research, and consider some consequences of each. Next, I provide critiques of some recurring assumptions about automated mathematics. Finally, I ask how automation might reshape economies of labor and credit in mathematics, and briefly offer my own hopes.

\section{Post-Cartesian and Post-Leibnizian Standards of Proof}

Mathematicians have never agreed on a single criterion that decisively verifies a claimed proof. Since the earliest aspirations to certainty in ancient times, the practice of mathematics has always involved the coexistence of multiple styles of proving \cite{netz_shaping_1999,netz_ludic_2009,chemla_history_2012}. A given demonstration can be deemed a “proof” by a group of people if it addresses all of their explicit and implicit epistemic fears, and different people worry about different things.

The history of modern mathematics is filled with failed attempts to convince all mathematicians to adopt a single and clear-cut criterion of proof. No such attempt has ever achieved consensus. Recent proposals for computer-based formal verification are only the latest episode in this history, and perhaps they will become the first successful exception. I suggest that past and present controversies about the status of computer-assisted proofs reflect a longstanding tension in modern mathematics, which may be traced to the distinction between two competing ideals of proof. Borrowing the terminology from philosopher Ian Hacking, we might call them the “Cartesian” and the “Leibnizian” ideals of proof.

\begin{description}
\item[Cartesian ideal of proof] “after some reflection and study, one totally understands the proof, and can get it in one’s mind ‘all at once’” \cite[p.~21]{hacking_why_2014}.
\item[Leibnizian ideal of proof] “every step is meticulously laid out, and can be checked, line by line, in a mechanical way” \cite[p.~21]{hacking_why_2014}.
\end{description}

These ideals are associated with René Descartes and Gottfried Leibniz, who were only two among many theorists of proof in the seventeenth century. Their contemporaries engaged in various debates about what should and should not count as a proof, including controversies about whether geometers should reject proofs by contradiction \cite{mancosu_philosophy_1996}. I focus on Cartesian and Leibnizian ideals because I find their contrast especially relevant to current discussions of automation.

Neither ideal has been achieved completely, much less required, in most actually existing mathematics. But mathematics has involved a mix of the two ideals---or, more precisely, a mix of practical standards that may be seen as derivative forms of these ideals. Both ideals repeatedly faced practical challenges, and had to be replaced by more realistic standards. To begin with, some proofs are too long and complex to fit in one’s mind “all at once.” In such cases, the Cartesian ideal had to give way to a less strict practical standard.

\begin{description}
\item[Post-Cartesian standard of long and complex proofs] after some study, one understands the proof’s general strategy and each of its parts, and can follow it one part at a time.
\end{description}

Another challenge is that some proofs are beyond the grasp of an isolated individual, and their construction and verification may require multiple people. This problem became more pressing with the professionalization and specialization of mathematics. After the nineteenth century, no individual mathematician could claim to know all branches of mathematics at the level of specialists, and sometimes Henri Poincaré is said to have been the last who could have aspired to do so. There are numerous famous examples of highly specialized proofs that demanded several years of efforts at explanation and persuasion \cite{dutilh_novaes_dialogical_2021}. These include ultimately accepted proofs, like that of Fermat’s Last Theorem, as well as controversial claimed proofs, like that of the $abc$ conjecture.

\begin{description}
\item[Post-Cartesian standard of highly specialized proofs] each part of the proof is understood, after some study, by some group of trusted specialists.
\end{description}

Proof verification is a sophisticated social process. A claimed proof often has parts that can be understood only by a few specialists in a subfield, even if its general strategy may be accessible to other professional mathematicians. In practice, the adjudication of specialized proofs relies on social institutions, such as journals, seminars, committees, degrees, awards, and symbolic markers of reputation and trust \cite{rosental_weaving_2008,greiffenhagen_materiality_2014,barany_myth_2015,harris_mathematics_2015}. In this regard, mathematics is not fundamentally different from other fields of knowledge. The production and verification of an experimental proof in physics---say, the detection of gravitational waves---may require hundreds of specialists not only in experimentation but also in theory and in instrumentation, all of whom must be trusted \cite{galison_image_1997,collins_gravitys_2017}.

The Leibnizian ideal has also developed into different practical standards. One key transformation was the turn to modern axiomatic proofs, which was advanced by David Hilbert’s axiomatization of geometry at the end of the nineteenth century, and later embraced by various mathematicians in the first decades of the twentieth \cite{corry_origins_1997,corry_david_2004}. The members of Bourbaki are closely associated with the turn to modern axiomatics \cite{corry_nicolas_1992,aubin_withering_1997}. As historian Alma Steingart has elegantly documented, modern axiomatics influenced broader intellectual transformations in the mid-twentieth century, including in structuralist social science and in abstract art \cite{steingart_axiomatics_2023}.

\begin{description}
\item[Post-Leibnizian standard of modern axiomatic proofs] every step can be checked, line by line, by means of strict logical deductions, and in principle can be derived from axiomatic foundations.
\end{description}

The key here is “in principle.” In practice, published proofs are hardly ever actually traced to the axioms. Moreover, what counts as the proper “foundations” has continually been the subject of dispute. Until the late nineteenth century, the foundations of geometry were thought to be grounded in physical reality, but this idea was challenged after the formulation of non-Euclidean geometry. In the early twentieth century, there were numerous controversies about the foundations of logic, particularly about the law of excluded middle, as well as efforts to devise alternative mathematical logics around the world, for example in the late Russian Empire \cite{tatarchenko_russian_2021}. Ever since Zermelo--Fraenkel set theory became a leading candidate for a unified foundation of mathematics, the inclusion of the axiom of choice has been controversial \cite{moore_zermelos_1982,mancosu_adventure_2010}. And as ongoing debates about univalent foundations and homotopy type theory make clear, there are still compelling reasons to keep questioning one’s choice of axioms \cite{the_univalent_foundations_program_homotopy_2013,pelayo_homotopy_2014,grayson_introduction_2018}.

Current practices of proof reflect a complex configuration of multiple standards, which remain largely implicit. Mathematicians learn these standards not by reading them explicitly in textbooks, but mainly through practice. I shall try to render them explicit by venturing a first approximation.

\begin{description}
\item[The current configuration of proof standards (approximation)] post-Leibnizian standards adjudicate the validity of a proof in principle, while post-Cartesian standards adjudicate validity in practice in addition to other considerations beyond validity.
\end{description}

Since the first half of the twentieth century, the post-Leibnizian standard of modern axiomatic proofs has come to define how mathematicians adjudicate the validity of a proof in principle. Even if a published proof may not include every detail of every step all the way from the axioms, mathematicians consider it valid if they believe that every step would be valid if checked. However, post-Cartesian standards remain a central part of mathematical practice, especially for long, complex, and highly specialized proofs. Few published proofs are traced explicitly to axiomatic foundations. In practice, mathematicians depend mainly on post-Cartesian standards to adjudicate the validity of a proof, involving social processes such as peer review, preprint circulation, and seminar discussion.

Post-Cartesian standards also play a role in additional considerations beyond validity. A proof is often said to be “better” or more “elegant” or more “interesting,” or to contribute to a deeper “understanding,” when it approximates the Cartesian ideal. Among equally valid proofs, mathematicians tend to prefer those that are shorter, can be understood by many people, and do not require computers.

\section{Forms of Automation in Mathematical Research}

\begin{figure}[H]
\centering
\vspace*{0.25cm}
\includestandalone[width=1.0\textwidth]{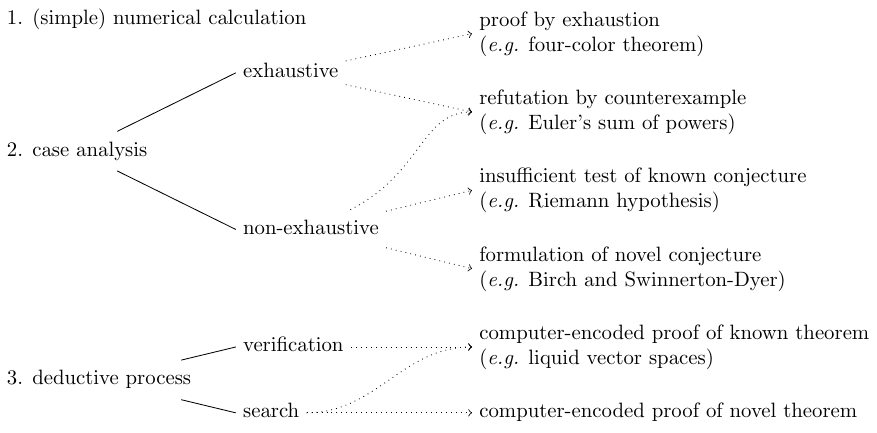}
\vspace*{0.25cm}
\caption{Diagram of prominent forms of automated mathematics.}
\label{diagram}
\end{figure}

How might computers transform current standards of proof? To address this question, it is helpful to distinguish between some forms of automation in mathematical research, involving different uses of computing technologies for different purposes. My intention is not to provide an exhaustive taxonomy, but simply to organize my discussion of prominent forms (see Figure \ref{diagram}).

I distinguish between three major categories. The first is the automation of simple \emph{numerical calculation}, the sort that can be done with an abacus, a slide rule, or a non-programmable electronic calculator. Insofar as this first category of automation is concerned, the emergence of digital computers in the mid-twentieth century constitutes an advance in degree rather than in kind: the ability to calculate larger numbers, but still one at a time. This sort of numerical calculation does not pose challenges to current standards of proof.

The second major category is the automation of \emph{case analysis}. Digital computers can be programmed to search large numbers of cases that humans alone cannot. Sometimes, a computer-assisted analysis can find a counterexample that disproves a conjecture. An early instance is the refutation of Euler's sum of powers conjecture by L. J. Lander and T. R. Parkin, who found a counterexample with the aid of a computer in 1966 \cite{lander_counterexample_1967}. A more recent instance is the refutation of Keller’s conjecture in 2019 \cite{brakensiek_resolution_2022}.

Other times, the analysis may not find any counterexample after searching all cases, and thus constitute a proof by exhaustion. The classic illustration is the proof of the four-color theorem. In 1976, Kenneth Appel and Wolfgang Haken reduced the original problem to less than two thousand cases, then checked each case with a computer \cite{appel_every_1989,mackenzie_mechanizing_2001}. A more recent illustration is the proof that every position of the Rubik’s Cube can be solved in 20 moves or fewer, published in 2010 \cite{rokicki_diameter_2014}. It involved reducing the roughly $4.3 \times 10^{19}$ positions into a smaller space of cases, which still took the equivalent of 35 years of the processing power of a desktop computer of the time to be searched. Computer-assisted proofs by exhaustion have triggered quite a bit of unease and resistance, because their acceptance demanded a new standard of proof that moved another step away from the Cartesian ideal of total understanding.

\begin{description}
\item[Practical standard of computer-assisted proofs by exhaustion] after some study, one understands the proof’s general strategy as well as the computer program that checked all the cases.
\end{description}

A case analysis may also be non-exhaustive, either because the number of cases is infinite or because it is finite but too large for existing computers. Some non-exhaustive analyses may serve mildly to strengthen confidence in a known conjecture. For instance, computers have tested trillions of nontrivial zeros of the Riemann zeta function $\zeta(s)$, all of which have real part $1/2$, in accordance with the Riemann hypothesis \cite{platt_riemann_2021}. Such non-exhaustive tests might be useful for certain practical applications, even if they are insufficient for mathematicians, who are well aware of instances in which the first exception to a conjecture occurs only after a fantastically large number of cases. For example, the earliest known case when $\pi(x) > \text{li}(x)$ (the prime-counting function exceeds the logarithmic integral function) happens around $x \sim 10^{316}$ \cite{saouter_sharp_2010}.

More compellingly, a non-exhaustive case analysis may also lead to the formulation of a novel conjecture. Instead of merely testing cases of a known conjecture, a researcher might run analyses in a more exploratory fashion, and discover unexpected patterns that hint toward a conjecture yet unknown. In the early 1960s, Bryan John Birch and Peter Swinnerton-Dyer used a computer to analyze some special cases of elliptic curves, leading to the formulation of their now famous conjecture \cite{birch_notes_1965}. In recent years, researchers have fruitfully deployed machine learning to formulate unprecedented conjectures in fields as varied as knot theory and representation theory \cite{davies_advancing_2021}. Such exploratory analyses leading to novel conjectures, which do not necessarily pose any challenge to current standards of proof, arguably constitute the least controversial and most compelling form of automation among those within the second category (case analysis).

Finally, the third major category is the automation of the \emph{deductive process} itself. This includes both the verification of known theorems and the search for novel ones. Starting in the late 1960s, software systems for formal verification appeared in different places: Automath, founded by N. G. de Bruijn in the Netherlands in 1967; Logic for Computable Functions (LCF), started by Robin Milner in California in 1969 and further developed in Scotland; and Mizar, founded by Andrzej Trybulec in Poland in 1973 \cite{nederpelt_selected_1994,paulson_logic_1987,grabowski_four_2015}. De Bruijn, a mathematician, has recalled his firsthand experience in checking a long and repetitive proof in combinatorics as a motivation for Automath. Meanwhile, Milner’s work on LCF emerged from his research program in theoretical computer science, and Trybulec’s initial vision for Mizar focused on a centrally managed database or “library” of mathematical papers. Many other projects followed over the decades, including Isabelle, HOL, Coq, Metamath, and Lean. Each of these projects was motivated at least in part by some version of the epistemic fear that conventional mathematical proofs, including highly trusted ones in the published literature, may contain flaws because their ever-increasing complexity and specialization makes them susceptible to human error. However, the exact philosophical motivations, axiomatic foundations, design choices, and narratives for each project have differed \cite{morgan_politics_2022}. Taken together, such projects have advanced a new standard of proof.

\begin{description}
\item[Practical standard of computer-encoded proofs] every step can be checked by a computer program and derived from the axiomatic foundations of the program; and after some study, one understands or trusts the encoding of the proven statement.
\end{description}

This standard has both post-Leibnizian and post-Cartesian dimensions. It is often assumed to be sufficient, because its post-Leibnizian dimension is seen as stricter than the conventional standard of axiomatic proofs in print, which is already accepted \cite{azzouni_tracking_2006,azzouni_why_2009,avigad_reliability_2021,hamami_mathematical_2022}. Some proponents of “machine-checkable” or “formal” proofs tend to emphasize this post-Leibnizian dimension, and may wish to eliminate post-Cartesian standards involving the elusive notion of “understanding.” I prefer to speak of computer-encoded proofs, both to avoid the ambiguity of the word “formal,” which once described axiomatic proofs which some people now call “informal,” and to clarify that machines are not the only agents involved in checking.

Computer-encoded proofs do not evade the need for human understanding. Humans must still read the code and verify that the proven statement, as encoded in the program, actually represents or corresponds to the claim under consideration. This issue might sound trivial, but can be quite complex and delicate. Consider the formalization of Peter Scholze’s definition of perfectoid spaces by Kevin Buzzard, Johan Commelin, and Patrick Massot in 2020 \cite{buzzard_formalising_2020}. Before the formalizers could state the definition in Lean code, they had to encode many other definitions and theorems in topology, algebra, and geometry, amounting to thousands of original lines of code. Buzzard, Commelin, and Massot worked hard to ascertain that the code was equivalent to the mathematical object that Scholze had in mind. The equivalence is far from obvious, and demanded close readings of both the code and Scholze’s original paper. Many researchers report that the exercise of encoding a result in systems like Lean, Coq, and Isabelle has improved their comprehension of it. The process is arduous and instructive.

Although computer-encoded proofs remain susceptible to human error and introduce new epistemic fears (such as software bugs), I have not yet encountered public objections to their assumed sufficiency as a standard of proof. A more controversial question is whether computer-encoded proofs should become not only sufficient but also necessary. Should mathematical journals require computer programs along with, or even instead of, the submitted papers? Should all theorems be demoted to conjectures until computer-encoded proofs are available? Although proposals of this kind have already appeared for several decades, they have garnered more attention in recent years. One reason is that computer-encoded proofs have become realistic for increasingly sophisticated results. In 2022, a team of more than twenty researchers succeeded in formalizing, in Lean, a particularly challenging theorem regarding liquid vector spaces by Dustin Clausen and Peter Scholze \cite{scholze_liquid_2022,castelvecchi_mathematicians_2021}. It is also worth noting that the difficulty of encoding theorems in a given system may vary depending on its axiomatic foundations. In 2021, Anthony Bordg, Lawrence Paulson, and Wenda Li formalized Grothendieck’s schemes in Isabelle, despite the difficulty of doing so with Isabelle’s higher-order logic, which is weaker than Lean’s dependent type theory \cite{bordg_simple_2022}.

In parallel with these efforts at verifying existing results, there are also attempts to develop more automated methods of searching for new proofs and theorems. The latter task is actually older. The earliest projects about automated deduction in the 1950s and 1960s had focused on search, not verification. The application of “artificial intelligence” to the search for proofs is not new either. Indeed, from its very beginnings, the field of AI has involved the aspiration to search for proofs automatically. At the 1956 Dartmouth Summer Research Project on Artificial Intelligence, where the term was coined, Allen Newell and Herbert Simon presented the Logic Theorist, a computer program that could prove theorems from Bertrand Russell and Alfred North Whitehead’s \emph{Principia Mathematica}, and proposed their influential view of human minds and digital computers as “species of the same genus,” namely symbolic information processing systems \cite{dick_artificial_2019,crowther-heyck_herbert_2005}. Early researchers of automated deduction found inspiration in diverse philosophical ideas. Historian Stephanie Dick has shown that Hao Wang, who also developed computer programs for proving theorems in the late 1950s, emphasized differences between humans and machines, as well as among human minds, due to his alternative theory of human reason inspired by Marxist philosophy \cite{dick_marxist_2023}.

Although some projects hoped to fully automate the process of theorem proving, they never really did. The work of “automated theorem proving” inevitably demands certain moments of human insight and intuition. This was already recognized by researchers in the 1970s and 1980s, as Dick has proven in her account of the Automated Reasoning Assistant (AURA) project at the Argonne National Laboratory \cite{dick_aftermath_2011}. In this sense, the boundary between “automated” and “interactive” theorem proving has always been blurry. Rather than two separate fields, I prefer to think of one field of partially automated practices involving multiple variations of human guidance and interactivity.

I am not aware of any important theorem that was first proven through a deeply automated search involving limited human guidance. However, it is possible that such a theorem may appear in the future. It is also possible that its derivation may be so inscrutable that no humans will be able to understand even its general structure. If this happens, I predict that a new candidate standard of proof will become the subject of intense debate.

\begin{description}
\item[Potential practical standard of deeply automated proofs] one trusts the computer programs that constructed and checked the proof, even if one may not understand it at all.
\end{description}

Will mathematicians accept this standard as sufficient? Only time can tell. Notice that this situation is different from that with computer-assisted proofs by exhaustion, where one understands at least the general strategy. It also differs from the situation with human-authored computer-encoded proofs, which are understood at least by their authors. I emphasize that the answer is not obvious or self-evident. Both post-Leibnizian and post-Cartesian standards play a role in mathematics, and sometimes they present competing considerations. To many mathematicians, understanding a proof or statement is not just a requirement for believing it, but the very point of doing mathematics.

\section{Critiques of Recurring Assumptions}

To recap, recent advances in the automation of deduction (the third category) present mathematicians with two key decisions about their standards of proof. One is whether computer-encoded proofs should become required. The other is whether completely inscrutable, deeply automated proofs should be accepted. I do not wish to defend a normative position about either question. Instead, I provide critiques of some assumptions about automated mathematics that recur in responses to these questions. My intention is to stimulate discussions that move beyond these assumptions, and to raise different questions.

The first assumption is that computer-encoded proofs can eliminate the need for social institutions of adjudication, such as peer review. Even in a scenario where computer-encoded proofs become required, they still demand some human involvement in the process of verification. As we have seen, at a minimum, they demand social consensus that a piece of code represents or corresponds to the claimed mathematical statement. Social institutions are still needed to adjudicate the equivalence of mathematical ideas in code, in print, and in mind.

The need for social institutions is even clearer in the adjudication of the importance and value of results. As Akshay Venkatesh has helpfully suggested, an open problem is more highly valued if it is considered more difficult and more central \cite{venkatesh_thoughts_2022}. It is difficult if many people try to solve it and fail; it is central if it is linked with many other questions of prior importance. Can the assessments of difficulty and centrality be made algorithmically? In my view, only to a limited extent. Perhaps the graph of dependencies among theorems in a proof library can give a rough indication of centrality, akin to a citation network for conventional papers. But citation counts and related metrics are notoriously poor indicators of the importance or value of mathematical papers. The social judgments of difficulty and centrality---not to mention other virtues like “beauty” and “originality”---can be nuanced in ways that elude automatic quantification.

Another recurring assumption is that it is feasible to formalize all mathematical knowledge. It is tempting to think of the task of formalization as merely one of “translating” the existing content of published papers into a formal language. But this characterization underestimates the difficulty of the task. To begin with, published papers do not usually contain every detail of every step. Moreover, much of mathematical knowledge consists of unpublished results that circulate informally among specialists, or “mathematical folklore” \cite{harel_folk_1980,rittberg_epistemic_2020}. Indeed, one of the potential benefits of formalization is to hold such folklore to account. In the course of their training, mathematicians gradually learn to judge which published or unpublished claims are trustworthy, as well as which details to include and which to exclude from a paper. This learning does not happen by following explicit or formal instructions from textbooks. Rather, it happens mainly through practice, over time. Much of mathematical knowledge is difficult to formalize or codify explicitly, amounting to what sociologists call “tacit knowledge” \cite{collins_tacit_2010}.

The boundaries of “mathematics” as a field of knowledge are not fixed or self-evident. Each formal language or software system prescribes its redefinition of those boundaries, which can never include all of mathematics. A certain version of this point is blatantly clear, since much of what is called “mathematics” outside academic settings and across the world, including “ethnomathematics,” does not fit the modern axiomatic form \cite{ascher_mathematics_2002}. But there are subtler versions of the point. For instance, diagrammatic proofs can be difficult, if not impossible, to express in a formal language \cite{azzouni_that_2013}. Moreover, the axiomatic foundations of a system determine the limits of what is expressible and provable within it. It is worth noting that the desire for computer-based formal verification has served to motivate the development and adoption of new foundations of mathematics. Vladimir Voevodsky acknowledged this motivation as central to his work on univalent foundations and homotopy type theory \cite{voevodsky_origins_2014}. All formal languages and software systems inevitably exclude some kinds of mathematical knowledge, and render some kinds more easily expressible and legible than others.

Lastly, I want to address the implicit assumption that the growth of automated mathematics is driven primarily by the intellectual needs or priorities of academic mathematicians. Rather, this area of research has been shaped significantly by the strategic interests of the military and corporate institutions that have funded it. As sociologist Donald MacKenzie has documented, early efforts in “formal verification” or “program proof” emerged largely out of practical concerns about computer security, with funding from the US military, particularly the National Security Agency (NSA) \cite{mackenzie_mechanizing_2001}. For instance, the military worried that a real computer program or operating system might not conform to its supposed mathematical model of “security.” By the early 1970s, the NSA perceived formal verification as central to computer security. Although military uses of computer-assisted proof have changed significantly in the fifty years since then, the US Air Force Office of Scientific Research (AFOSR) still sponsors research in this area under its “Information Assurance and Cybersecurity” program, mentioning “interactive and automated theorem proving” as well as “formalized mathematics” among research areas that “are expected to provide valuable insights into various cybersecurity problems” \cite{us_air_force_research_laboratory_information_2020}.

In recent years, commercial firms such as Microsoft, Google, and Meta have become growing sources of investment in automated mathematics. Unlike the military, which focuses on specific problems with more predictable applications, those firms benefit from their investments in less direct ways. One way is to attract mathematically talented researchers to corporate-owned platforms and tools for software development. Consider the case of Lean, which is developed primarily by Microsoft \cite{felty_lean_2015}. Although the beloved proof assistant is open-source software, its default interface is implemented as an extension to Microsoft’s code editor, Visual Studio Code, which is integrated tightly with other Microsoft products, such as “cloud computing services,” “integrated development environments,” and “software development kits” for Microsoft platforms and operating systems.

Commercial firms also invest in the automation of mathematics as part of a broader strategy of demonstrating AI capabilities in general, thus attracting public attention and capital investment. If the military tends to sponsor work on computer-assisted proofs under the rubric of “cybersecurity,” those firms have branded such work more often as “AI,” focusing lately on projects involving “deep learning” and “large language models.” Unsurprisingly, public-relations teams and corporate-sponsored researchers deploy hyperbolic language to publicize their work, for example by framing their research as an effort to “create an automatic mathematician” (Google) or even asking sensationally, “Is human led mathematics over?” (Meta) \cite{platzer_towards_2021,meta_ai_teaching_2022}. It remains to be seen whether academic mathematicians will embrace or refuse the modus operandi of hype-driven stock markets and venture capital.

\section{Economies of Labor and Credit}

Despite popular narratives of technological unemployment, automation does not simply “replace” human labor \cite{benanav_automation_2020}. Rather, automation displaces and reconfigures labor. It renders some kinds of work obsolete, creates some new kinds, and transforms the value of both the old and the new kinds. Anthropologists, sociologists, historians, and other scholars have studied how automation can reinforce gendered, racial, and transnational divisions of labor, and how seemingly automated technologies often rely on invisible, underpaid, and uncredited work, such as generating and cleaning the data used to train and validate machine-learning models \cite{nakamura_indigenous_2014,irani_cultural_2015,amrute_encoding_2016,hicks_programmed_2018,roberts_behind_2019,atanasoski_surrogate_2019}. In this final section, I ask how automation might reshape the economies of labor and credit in mathematics.

The economy of credit in mathematical research typically follows a first-past-the-finish-line criterion, at least for theorems. The first person or team to publish a document that includes the final steps to arrive at a theorem, satisfying current standards of proof, receives the central credit for the result. Perhaps this system does not do justice to the collective and gradual character of research, since every major result builds upon countless prior results, techniques, and ideas. The practices of attributing credit for the formulation of conjectures or research programs (like the Langlands program) are less consistent. An illustrative dispute concerned the attribution of the conjecture on the modularity of elliptic curves (now a theorem). The dispute reflected conflicting views regarding which contribution deserves priority: the first expression of an expectation along the lines of the conjecture (even if imprecise), the first statement of a version of the conjecture in private communication, the first publication of a precise statement, or the first elaboration of certain forms of evidence. In the absence of a consistent system of credit, some attitudes to this dispute seem to have been inflected by interpersonal conflicts and even national allegiances \cite{harris_virtues_2020}.

There is no consistent system for rewarding the work of producing a computer-encoded proof of an existing theorem. What might happen to the economy of credit if computer-encoded proofs become required in the future? One possibility is that the credit system will follow suit: that the first people to provide a computer-encoded proof would receive the central credit for a theorem. It is unlikely that such a change would apply retroactively, because this has not usually happened with past changes. We still say “Euler’s formula” even though his original demonstration would not satisfy later standards of proof, which have changed many times since the eighteenth century. Another possibility is that the economy of credit and the standards of proof will become misaligned: the central credit would still go to the author of the first print publication containing a demonstration, but the result would only become accepted after someone (not necessarily the same person) provides a computer-encoded proof.

An even harder challenge posed by automation is the emergence of large collaborations involving multiple roles which transcend the traditional expertise of mathematicians. To be sure, mathematicians have long engaged in large collaborations among themselves. Decades before the rise of “massively collaborative mathematics” on the internet, the classification of finite simple groups involved more than one hundred mathematicians, who authored three to five hundred articles until the proof was declared completed in 1981 \cite{gowers_massively_2009,steingart_group_2012}. The involvement of non-mathematician co-authors in mathematical research is not new either. Applied mathematicians already need to negotiate authorship with collaborators from other disciplines, like biologists and economists, who may provide the motivation, inspiration, or application for a mathematical result but who might not contribute directly to the strictly mathematical work. What is distinct about automation, particularly in the search for new proofs and theorems, is the proliferation of roles---such as programming computers, optimizing algorithms, tweaking parameters, and cleaning data---which not only transcend the traditional expertise of mathematicians but also contribute to the very process of mathematical reasoning, even in pure mathematics.

This proliferation of roles will raise difficult questions regarding collective authorship, similar to those that the experimental sciences have already faced for several decades. To give an illustration from particle physics, the planners of a new detector at the Stanford Linear Accelerator Center (SLAC) devised authorship protocols in 1988 which not only varied depending on the type of publication (internal memos, conference proceedings, review papers) but also established a hierarchical structure: “physicists” were included as authors in all “physics papers,” while “engineers” were excluded from most \cite{biagioli_collective_2003}. It is plausible, but perhaps undesirable, that mathematical collaborations may reproduce analogous hierarchies reinforcing disciplinary boundaries, for example between “mathematicians,” “computer scientists,” “software engineers,” and “research assistants.” In biomedicine, questions of authorship are tangled with corporate battles over intellectual property, including profitable pharmaceutical patents \cite{mirowski_instability_2002,sunder_rajan_pharmocracy_2017}. Hopefully, mathematics will be able to avoid some of these issues, since theorems and conjectures are not ordinarily patentable. But the field will still need to grapple with the politics of attribution, compensation, and professional recognition.

If the process of searching for new theorems and proofs becomes more automated, there will be an increasing need for interpreting computer-generated proofs and explaining them to people. This work might be aided by the future development of techniques that facilitate human interpretation and explanation, for example by visualizing the general structure of proofs. However, the existing techniques of “interpretable” and “explainable” AI and machine learning, such as those promoted by the US Defense Advanced Research Projects Agency (DARPA) through its Explainable Artificial Intelligence Program, do not seem to address the kinds of interpretability and explainability that mathematicians would need \cite{rudin_interpretable_2022,gunning_darpas_2019}. Many existing techniques aim to translate probabilistic and statistical models into a sequence of logical deductions, such as decision trees that produce approximately the same outputs given the same inputs. But in the case of automated proofs, the uninterpretable object is already in a deductive form. The problem is that the sequence of deductions is itself difficult to interpret or explain. Until specialized solutions to this problem exist, I expect that human labor will remain vital. Whether this explanatory labor will be properly valued, when decoupled from the credit for conjecturing or proving the explicandum, I cannot say.

I also worry that raw computational power might become a significant source of competitive advantage in searching for theorems and proofs. If so, the reliance on automated methods may amplify the advantage of wealthier research institutions, as well as researchers’ dependence on corporate and military sponsorship for access to supercomputers. To counter this risk, mathematicians may need to reform their economy of credit. In particular, they may need to shift emphasis away from the first-past-the-finish-line criterion for theorem credits, and toward the recognition of the gradual and collective labor that already underlies all major results.

Thus, mathematical researchers face many collective decisions and possible futures. It is quite plausible that automation may lead to unjust labor hierarchies, reinforced disciplinary boundaries, and greater dependence on military and corporate sponsorship. But none of these plausible consequences are inevitable. Researchers, at least in academic settings, can choose to shape their own future otherwise. This will demand a conscious and organized struggle against some of the external forces and powerful interests at play. My hope is that researchers will reap the benefits of new technologies cautiously and selectively---while working actively to resist corporate and military influence, to dismantle labor hierarchies, to communicate across disciplinary boundaries, to implement more equitable systems of credit and compensation for all contributions, and to celebrate the plurality of available ways of knowing and doing mathematics.

\bibliographystyle{amsplain}

\bibliography{main}

\end{document}